\documentclass{article}

\title{The generic special scroll of genus $g$ in $\P^N$. Special Scrolls in $\P^3$.}
\author{Luis Fuentes Garc\'{\i}a
\and Manuel Pedreira P\'erez\thanks{Partially supported by Xunta de Galicia, project PGIDITOPXIA20702PR}}
\date{}

\newtheorem{teo}{Theorem}[section]

\newtheorem{prop}[teo]{Proposition}
\newtheorem{cor}[teo]{Corollary}

\newtheorem{rem}[teo]{Remark}

\def\cosa{{}}

\def\ov{\overline}

\def\E{{\cal E}}
\def\Te{{\cal O}}
\def\R{{\bf R}}

\def\F{{\cal F}}

\font\euf=eufm10 at 12pt
\def\e{\mbox{\euf e}}
\def\b{\mbox{\euf b}}
\def\c{\mbox{\euf c}}
\def\aa{\mbox{\euf a}}

\def\P{{\bf P}}

\def\k{{\cal K}}

\def\qed{\hspace{\fill}$\rule{2mm}{2mm}$}

\def\lrw{{\longrightarrow}}

\begin{document}

\maketitle

{\footnotesize{\bf Authors' address:} Departamento de M\'etodos Matem\'aticos y Representaci\'on, Universidad de A Coru\~na. $15192$ Coru\~na. Spain. {\tt lfuentes@udc.es}\\ Departamento de Algebra, Universidad de Santiago
de Compostela. $15706$ Santiago de Compostela. Spain. e-mail: {\tt
mpedreira@usc.es} \\
{\bf Abstract:} We study the generic linearly normal special scroll of genus $g$ in $\P^N$. Moreover, we
give a complete classification of the linearly normal scrolls in $\P^3$ of genus $2$ and $3$. 
}

\vspace{0.1cm}

\section{Introduction.}
 Let $\pi:S=\P(\E)\lrw X$ be a geometrically ruled surface; consider a base
point free unisecant linear system $|H|$ on $S$. The image of $S$ by the regular map defined by the
complete linear system $|H|$ is a linearly normal scroll $R\subset \P^{h^0(\Te_S(H))-1}$. We define the
speciality of $R$ as $i=h^1(\Te_S(H))$. The scroll is special when $i>0$.

In this paper we study the linearly normal special scrolls. First we describe the generic linearly normal special scroll in $P^N$. The Hilbert Scheme of the nonspecial ruled surfaces was studied by the second author. He proved that there is an unique irreducible component corresponding to the smooth nonspecial linearly normal scrolls of genus g and degree $d$ if and only if $d\geq 2g+2$, this component dominates the moduli space of smooth curves $M_g$, and the generic nonspecial scroll is of maximal rank. These results are contained in \cite{Pe1} and \cite{Pe2}. Moreover, it was proved in \cite{Pe1} that the generic nonspecial scroll in $\P^N$ corresponds to stable geometrically ruled surfaces with invariant $e=-g$
or $e=-(g-1)$, depending of the parity of $g$ and $d$. These are the general type ruled surfaces studied by Segre in \cite{segre} and Ghione in \cite{ghione2}.

However, in general, the generic linearly normal special scroll $R$ of degree $d$ and genus $g$ in $\P^N$
does not correspond to the ruled surface of general type. In fact, by the Clifford Theorem for
stable vector bundles if $d\geq 4g-4$ the special ruled surface can not be stable. 

We will give a very detailed description of the ruled surface corresponding to the generic special
scroll. In particular, we see that it corresponds to a ruled surface of
general type precisely when $d\leq 3g-3$. Our method will be to project an adequate generic
decomposable special scroll.

The second part of the paper is pointed to study the linearly normal special scrolls in $\P^3$ of low
genus. Edge in \cite{edge} give a classification by studying curves in the
Grassmannian of lines of $\P^3$. He only considers curves of degree less than $7$. We give the complete
classification of genus $2$ and $3$, combining his methods with the study of the projections of
suitable decomposable special scrolls.

\bigskip

\section{Notation.}

We will follow the notation of \cite{hartshorne}.  We refer to \cite{fuentes} for a systematic
development of the projective theory of scrolls and ruled surfaces that we will use in this paper. In
particular:

Given a geometrically ruled surface $\pi:S\lrw X$ we know that $S=\P(\E_0)$, where $\E_0$ is the normalized locally free sheaf of rank $2$ over $X$ (see \cite{hartshorne}, V.2). We denote by $\e$ the invertible sheaf $\bigwedge^2\E_0$ and $e=-deg(\e)$. The geometrically ruled surface $S$ is called decomposable if $\E_0=\P(\Te_X\oplus \Te_X(\e))$ and it is called indecomposable in other case. 

A linearly normal scroll $R$ in $\P^N$ is given by a geometrically ruled surface $S$ and a
base point free linear system $|H|$ in $S$. In this way, $R$ is the image of $S$ by the map
defined by $|H|$. We say that the pair $(S,|H|)$ corresponds to the scroll $R$. If $x$ and $C$ are respectively a point and curve in $S$, the image of $x$ and $C$ in $R$ will be denoted by $\ov{x}$ and $\ov{C}$.

We will denote the canonical divisor of a smooth curve $X$ by $\k$. Given a divisor $\b\in Div(X)$ of degree $\geq 2g-2$, we denote the decomposable ruled surface
$\P(\Te_X\oplus
\Te_X(\b-\k))$ by
$S_{\b}$ and the scroll defined by the linear system $|X_0+\b f|$ by $R_{\b}$.

If $\aa$ is a divisor in $X$ and $C$ a curve in a ruled surface $S$ we denote by $C+|\aa|f$
the set of reducible curves $\{ C+\aa' f/ \aa\sim \aa'\}$.

Given two subsets $A,B$ of a linear system $|H|$, we denote by $\langle A,B \rangle$ the minimum
linear subsystem in $|H|$ containing $A$ and $B$.

\section{Special scrolls in $\P^N$.}\label{especialesN}

Let $R\subset \P^N$ be a linearly normal special scroll of degree $d$, genus $g$ and
speciality $i$. By the Riemann-Roch Theorem:
$$
d-2g+1+i=N.
$$
Because $R$ is special, $i\geq 1$; so we deduce that:
$$
d\leq 2g+N-2.
$$

\begin{prop}\label{proyeccion}
Any linearly normal special scroll $R\subset \P^N$ of degree $d$, genus $g$ is
the projection of a decomposable scroll $R_{\b}\subset P^{2g+N-3}$ from
$4g+N-5-d$ points spanning a $\P^{2g-4}$; $\b$ is a nonspecial divisor of degree $2g+N-3$.
\end{prop}
{\bf Proof:}
Because $R$ is special it has a special curve $\ov{X_a}$ (see \cite{pedreira2}). This is the
projection of a canonical curve. On the other hand, $d\leq 2g+N-2$, so  there
exist a unisecant irreducible curve $\ov{X_b}$ of degree $d-1\leq 2g+N-3$. It can be
projected from a nonspecial curve of degree $2g+N-3$. From this the scroll $R$ is the
projection of the scroll
$R_{\b}$. \qed

\begin{rem}

We can obtain the scroll $R_{\b}$ explicitly. Suppose that $R$ is given by the pair $(S,|H|)$.
Let $\Te_{X_a}(H)\cong \Te_X(\aa_1)$, with $\aa_1$ a special divisor and
$\Te_{X_b}(H)\cong\Te_X(\aa_2)$, with $deg(\aa_2)\leq 2g+N-3$. Let
$P_1+\ldots+P_l\sim \k-\aa_1$ and $\b\sim Q_1+\ldots+Q_m+\aa_2$ a nonspecial and birational
very ample divisor of degree $2g+N-3$. We make the elementary
transformation of $S$ at the points of $X_a\cap X_b$, $X_a\cap (Q_1+\ldots+Q_m)f$ and
$X_b\cap (P_1+\ldots+P_l)f$. We have:
$$
\begin{array}{l}
{X_a'.X_b'=X_a.X_b-\#(X_a\cap X_b)=0,}\\
{\pi_*(X_a'\cap H')\sim \pi_*(X_a\cap H)+P_1+\ldots+P_l\sim \k,}\\
{\pi_*(X_b'\cap H')\sim \pi_*(X_b\cap H)+Q_1+\ldots+Q_m\sim \b,}\\
\end{array}
$$
where $H$ is a generic hyperplane section. We obtain the decomposable ruled surface
$S_{\b}\cong
\P(\Te_X(\k)\oplus \Te_X(\b))$. The scroll $R_{\b}$ is the image of $S_{\b}$ by the map defined
by the linear system $|H'|$.
\end{rem}

\begin{rem}

Remark that, in general, $R_{\b}$ is not a canonical scroll in the sense of \cite{pedreira2}. Moreover, if $R_{\b}$ is a
canonical scroll the points of projection could not be in the canonical bisecant curve. 

\end{rem}

As a consequence of this Proposition, the generic special scroll in $\P^N$ will be obtained by
projecting a generic $R_{\b}$ from $k$ points in general position.  Since they span a
$\P^{2g-4}$, $k=2g-3$ and the projection has degree $d=N+2g-2$ and speciality $1$.

Note that to project $R_{\b}$ is equivalent to make the elementary transform of
$S_{\b}$. The invariant $-\e$ of this ruled surface is a generic
divisor of degree $e=deg(\b-\k)$. By (\cite{fuentes}, Theorem $50$) if
$e\geq g$ then the elementary transformation of $S_{\b}$ at a generic point $x\in Pf$ is
$S_{\b-P}$. From this, we deduce the following results:

\begin{prop}\label{general}
If $N\geq 3g-3$ then the generic linearly normal special scroll of genus $g$ and speciality
$1$ in $\P^N$ is a decomposable scroll $R_{\b}$, where $\b$ is a generic divisor of degree $N$.
\end{prop}

\begin{prop}\label{projection}
If $N<3g-3$ then the generic linearly normal special scroll of genus $g$ and speciality
$1$ in $\P^N$ is the projection of a decomposable scroll $R_{\b}$, where $\b$ is a generic
divisor of degree $3g-3$.
\end{prop}

Now, we have to study the projections of the decomposable scroll $R_{\b}$, where $\b$ is a generic
divisor of degree $3g-3$. $R_{\b}$ be the image of $S_{\b}$ by the map
defined by the linear system $|H|=|X_0+\b f|$. The invariant $\e$ of $S_{\b}$ is a divisor $\e\sim
\k-\b$ verifying $h^0(\Te_X(-\e))=0$.

The ruled surface $S_{\b}$ has the
following families of unisecant curves with generic irreducible element (see \cite{fuentes}):
\begin{enumerate}

\item $X_0$, with $X_0^2=1-g$ and $dim(|X_0|)=0$.

\item $X_1$, with $X_1^2=g-1$ and $dim(|X_1|)=h^0(\Te_X(-\e))=0$.

\item Given $c>0$, we call ${\cal F}_c=\{Y_c\subset S_{\b}|\,Y_c\equiv X_1+cf \}$, where $Y_c^2=g-1+2c$.

\end{enumerate}

\begin{prop}

The generic curve $Y_c\in {\cal F}_c$ is irreducible. Moreover, $dim{\cal F}_c=2c$.

\end{prop}
{\bf Proof:} We define the following map:
$$
p:{\cal F}_c\lrw Sym^c(X); \qquad p(Y_c)=\pi^*(Y_c\cap X_0).
$$
The fiber over any set of points $\c=P_1+\ldots+P_c\in Sym^c(X)$ is:
$$
p^{-1}(\c)=\langle X_1+\c f, X_0+|-\e+\c| f\rangle.
$$
with, $dim(p^{-1}(\c))=h^0(\Te_X(-\e-\c))=c$. We deduce that:
$$
dim{\cal F}_c=dim(Sym_c(X))+dim(p^{-1}(\c))=2c.
$$
Finally, the divisors $\b$ and $\k$ are effective without common base points. Applying the Theorem 30 of \cite{fuentes}, we obtain the irreducibility of the generic curve of ${\cal F}_c$.  \qed

\begin{cor}\label{conditions}
Let $\{x_1,\ldots,x_k\}$ be $k$ generic points in $S_{\b}$. If $2c\leq k$ then there is
a family of dimension $2c-k$ of irreducible curves in $\F_{c}$ passing through these points.
If $2c<k$ then there is not any irreducible curve in $\F_{c}$ passing through these points. 
\end{cor}

\begin{cor}\label{minimum}

Let $Z=\{x_1,\ldots,x_k\}$ be a set of $k$ generic points in $S_{\b}$. Suppose that $k=2c-j$, with $j\in \{0,1\}$. Let $S_{\b}'$ the elementary transform of $S_{\b}$ at $Z$. Then:
\begin{enumerate}

\item If $k\leq 2g-2$ then the minimum self-intersection curve of $S_{\b}'$ is $X_0'$, with
$X_0'^2=1-g+k$. There is a $j$-dimensional family of irreducible curves with self-intersection
$g-1+j$. Any other irreducible curve of $S_{\b}'$ has self-intersection greater than $g-1+j$.

\item If $k> 2g-2$ then there is a $j$-dimensional family of irreducible curves with
minimum self-intersection $g-1+j$. Any other irreducible curve of $S_{\b}'$ has
self-intersection greater than $g-1+j$.

\end{enumerate}

In particular, $S_{\b}'$ is not decomposable.

\end{cor}
{\bf Proof:} We use the properties of the elementary transformation (see \cite{fuentes}, Proposition 42). Since the points of $Z$ are generic:
$$X_0'^2=X_0^2+k=1-g+k; \qquad X_1'^2=X_1^2+k=g-1+k.$$
Let $Y_{c}\in \F_c$ be  a irreducible curve of $S_{\b}$. We have seen that $Y_{c}^2=g-1+2c$.
Suppose that $\{x_1,\ldots,x_l\}\subset Y_{c}$, with $l\leq k$. Then,
$Y_{c}'^2=Y_{c}^2-l+(k-l)=g-1+2(c-l)+k$. Moreover, by the Corollary \ref{conditions} we know that:

If $k\leq 2c$ then there is a family of dimension $2c-k$ of irreducible curves in $\F_c$
passing through $\{x_1,\ldots,x_k\}$. Thus, $Y_{c}'^2> g-1+2c-k$ for the generic
curve and $Y_{c}'^2=g-1+2c-k$ for a $(2c-k)$-dimensional family.

If $k>2c$ then $Y_{c}'^2=g-1-2c+k$ for a finite number of irreducible curves in $\F_c$ and
$Y_{c}'^2>g+1-2c+k$ for the generic curve.

Comparing the self-intersection of all these curves we complete the proof. Note that a ruled surface is decomposable if and only if there are two irreducible unisecant curves with the sum of their self-intersections equal to $0$. Therefore, in these cases the ruled surface $S_{\b}'$ is not decomposable. \qed

From this we conclude the following theorem:

\begin{teo}\label{generica}
Let $R\subset P^N$ be the generic special scroll over a curve $X$ of genus $g$. $R$ has
speciality $1$ and it has an unique special directrix curve. This curve is a canonical
curve. It is linearly normal if and only if $N\geq g-1$, equivalently when it is the curve of minimum degree of the scroll. In particular,

\begin{enumerate}

\item If $N\geq 3g-3$ ($d\geq 5g-5$) then $R$ is the image of the decomposable ruled surface
$\P(\Te_X\oplus
\Te_X(\k-\b))$ by the map defined by the linear system $|X_0+\b f|$, where $\b$ is a
nonspecial generic divisor of degree $N$. The special curve of $R$ is the image of the curve of minimum self intersection $X_0$.

\item If $g-1\leq N<3g-3$ ($3g-3\leq d< 5g-5$)then $R$ is the image of an indecomposable ruled surface
$\P(\E_0)$ with invariant $e=-(2g-2-N)$ and $\Te_X(\e)\cong \bigwedge ^2\E_0$, by the map defined by the
linear system $|X_0+(\k-\e)f|$. The special curve of $R$ is the image of the curve of minimum self intersection $X_0$.

\item If $3\leq N<g-1$ ($2g+1\leq d< 3g-3$) then $R$ is the image of an indecomposable ruled surface
$\P(\E_0)$ with $\Te_X(\e)\cong \bigwedge ^2\E_0$, by the map defined by the linear system
$|X_0+(\k-\aa)f|$, where $\aa\sim \pi_*(X_0\cap Y)$ and $Y$ is a unisecant irreducible curve of
self-intersection $2g-2+N$. The image of $Y$ in $R$ is the unique special directrix curve in
the scroll. It is a non linearly normal canonical curve. Moreover, 

\begin{enumerate}

\item if $N\equiv_2 g-1$ then $e=-(g-1)$.

\item if $N\equiv_2 g$ then $e=-g$.

\end{enumerate}

\end{enumerate}

\end{teo}
{\bf Proof:} If $N\geq 3g-3$, the conclusion follows from the Proposition \ref{general}. If $N<3g-3$, by the propositions \ref{proyeccion} and \ref{projection} we know that $R$ is the projection of $R_{\b}$ from $k=3g-3-N$ general points spanning a $\P^{k-1}$. Since, the unique special curve of $R_{\b}$ is a linearly normal canonical curve, the unique special curve of $R$ will be its projection. Moreover, this special curve is linearly normal if and only if:
$$
(k-1)+(g-1)<3g-3 \iff N\geq g-1.
$$
On the other hand, by the Corollary \ref{minimum} we know that $R$ is an indecomposable scroll. Furthermore, we can identify the curve of minimun degree of $R$.

If $k\leq 2g-2$ $(N\geq g-1)$, then the curve $\overline{X_0}$ of minimun degree of $R$ is the projection of the special curve of $R_{\b}$. In particular, $X_0^2=1-g+k=2g-2-N$.

If $k> 2g-2$ $(N< g-1)$, then the curve $\overline{X_0}$ of minimun degree of $R$ is not the projection of the special curve of $R_{\b}$. Now, $X_0^2=g-1+j$ where $j\equiv_2 k$, $j\in \{0,1\}$. Equivalently,
\begin{enumerate}

\item if $N\equiv_2 g-1$ then $X_0^2=(g-1)$.

\item if $N\equiv_2 g$ then $X_0^2=g$.

\end{enumerate} \qed

\section{Special scrolls in $\P^3$}\label{especiales3}

\subsection{General facts.}

Let $R\subset \P^3$ be a linearly normal special scroll of degree $d$, genus $g$ and
speciality $i$. Note that
$$
d-2g+1+i=3
$$
Because $R$ is special, $i>0$, so we deduce that:
$$
d<2g+2
$$
We apply the Proposition \ref{proyeccion} to the case $N=3$. We obtain:
\begin{prop}\label{proyeccion3}
Any linearly normal special scroll $R\subset \P^3$ of degree $d$, genus $g$ is the projection of a
decomposable the scroll $R_{\b}\subset \P^{2g}$ from $4g-2-d$ points spanning a $\P^{2g-4}$; 
$\b$ is a nonspecial divisor of degree $2g$.
\end{prop}

Let us study the cases of genus $2$ and $3$.

\subsection{Special scrolls of genus $2$}\label{genero2}

We apply the Proposition \ref{proyeccion3} to the case of genus $g=2$:

\begin{prop}\label{proyeccion32}
Any linearly normal special ruled surface $R\subset \P^3$ of genus $2$ and degree $d$ is the
projection from a decomposable scroll $R_{\b}\subset P^4$ from a point of multiplicity $6-d$; $\b$ is
a divisor of degree $4$.
\end{prop}

Let us give a detailed description of the scroll $R_{\b}$ of genus $2$ and $deg(\b)=4$. It is the image of the decomposable scroll $S_{\b}$ by the map $\phi$ defined by the linear system $|X_0+\b f|$.

 If $\b=\k$ then $\phi$ is not birational. In this case $R_{\b}$ is a rational ruled surface. So we can discard this case and we will take $\b=\k+A_1+A_2$ with $A_1+A_2\neq \k$. The singular locus of $R_{\b}$ is a line $\overline{X_0}$ and an isolated singular point (the image by $\phi$ of the base point of $|X_1|$). 

Now, we study the projection of the scroll $p:R_{\b}\lrw R$ from a point $\ov{x}$. First let us suppose that $\ov{x}$ is singular point of $R_{\b}$:

\begin{enumerate}

\item If $\ov{x}\in \overline{X_0}$, then the line $X_0$ projects into a point. Then we obtain a cone over a
nonspecial plane curve of degree $4$ and genus $2$.

\item If $\ov{x}$ is the singular point of $\overline{X_1}$ then this curve projects into a double line. In fact, the scroll $R_{\b}$ projects into a rational ruled surface. The projection map is not birational, so we does not consider this case.

\end{enumerate}

If $\ov{x}$ is a nonsingular point the study of the projection map $p$ is equivalent to study the elementary transform of $S_{\b}$ at $x$:
$$
\begin{array}{ccc}
{S_{\b}}&{\lrw}&{R_{\b}}\\
{\downarrow}&{}&{\downarrow}\\
{S'_{\b}}&{\lrw}&{R}\\
\end{array}
$$
Applying the Theorem 50 of \cite{fuentes}, we obtain:

\begin{enumerate}

\item If $x\not\in X_0$ and $x\not\in A_1f\cup A_2f$. Then the ruled surface $S'_{\b}$ is $\Te_X\oplus \Te_X(\k-\b-P)$.

\item If $x\not\in X_0\cup X_1$ and $x\in A_1f\cup A_2f$. Then we obtain an indecomposable scroll in $\P^3$ with invariant $\e'=\k-P$, where $P=\pi(x)$.

\end{enumerate}

From this we obtain the following models of linearly normal special scrolls in $\P^3$ of genus $2$.

\begin{teo}

Let $R\subset \P^3$ be a linearly normal special scroll of genus $2$, degree $d$ and speciality $i$.
Then $R$ is isomorphic to one of the following models:

\footnotesize
$$
\begin{array}{|c|c|c|c|c|c|}
\hline
{d}&{i}&{\E_0}&{e}&{|H|}&{
\begin{array}{c}
{\hbox{Curves of minimun degree}}\\
{nC_k^*; \hbox{$k=$degree; $^*$=special}}\\
{\hbox{$n=$multiplicity of $C_k$ in $R$}}\\
\end{array}}\\
\hline
\hline
{4}&{2}&{\begin{array}{l}
{\Te_X\oplus\Te_X(-\b)}\\
{\b\not\sim 2\k}\\
{\b\sim \k+P+Q}\\
\end{array}}&{4}&{
\begin{array}{c}
{|X_0+\b f|}\\
{\hbox{(cone)}}\\
\end{array}}&{C_4\subset \P^2\quad (\infty^3)}\\
\hline
\hline
{5}&{1}&{\begin{array}{l}
{\hbox{indecomposable}}\\
{\e=-P}\\
\end{array}}
&{1}&{|X_0+(\k+P)f|}&{
\begin{array}{l}
{\phi(X_0)=2C_1^*\subset \P^1}\\
{C_4\subset \P^2\quad (\infty^2)}\\
\end{array}}\\
\hline
{5}&{1}&{\begin{array}{l}
{\Te_X\oplus\Te_X(\k-\b)}\\
{deg(\b)=3}\\
{\b\hbox{ base-point-free}}\\
\end{array}}&{1}&{|X_0+\b f|}&{
\begin{array}{l}
{\phi(X_0)=2C_1^*\subset \P^1}\\
{\phi(X_1)=3C_1^*\subset \P^1}\\
\end{array}}\\
\hline
\end{array}
$$
\normalsize
where $\phi$ denotes the map defined by the linear system $|H|$:
$$
\phi:\P(\E_0)\lrw R\subset \P^3.
$$

\end{teo}

\bigskip

\subsection{Special scrolls of genus $3$.}\label{genero3}

We apply the Proposition \ref{proyeccion3} to the case of genus $g=3$.

\begin{prop}\label{proyeccion33}
Any linearly normal special ruled surface $R\subset \P^3$ of genus $3$ and degree $d$ is the
projection from a decomposable scroll $R_{\b}\subset P^6$ from $10-d$ points laying in a plane;
$\b$ is a divisor of degree $6$.
\end{prop}

If $d\leq 6$ we have to describe the projection of $R_{\b}$ from $10-d\geq 4$ points. This is not easy if the points are not generic. To solve this problem, we will use another method.  We will study the curve in the Grassmannian $G(1,3)$ parameterizing the special scroll $\R^3$ (see \cite{edge}).

\subsection{Curves in the Grassmannian parameterizing scrolls in $\P^3$.}

It is well known that any scroll in $\P^3$ corresponds to a curve in the Grassmannian
$G(1,3)$ of lines of
$\P^3$. Let us recall some basic facts about this variety (see \cite{GiSo98}). 

It can be realized like a
smooth quartic $Q_2\subset
\P^5$.

It has two $3$-dimensional families of planes: the
$\alpha$-planes, corresponding to the lines passing through a point; the $\beta$-planes
corresponding to the lines contained in a plane.

The intersection of $\P^3$ with the Grassmannian is a smooth quadric, a quadric cone or two
planes. In the first case, the points of the quadric correspond to the lines joining two
disjoint lines of $\P^3$. In the second case, the points of the cone correspond to lines
meeting the line defined by the vertex of the cone.

Now, let $C$ be a curve of genus $g$ and degree $d$ in the Grassmannian. It defines a scroll if
$\P^3$ with the same genus and degree. If the curve is linearly normal then the scroll is
linearly normal. The reciprocal statement is not necessary true. Let us suppose that
$g>0$; then:

\begin{enumerate}

\item If $C\subset \P^2$ then $C$ lays either in a $\alpha$-plane or in a $\beta$-plane. In the
first case, it corresponds to a cone over the curve $C$. The cone is linearly normal if and
only if the curve $C$ is linearly normal. In the second case, the scroll lay in a plane, so it
is a degenerate case.

\item If $C\subset \P^3$ then the curve lies either in a smooth quadric surface or in a cone.

In the first case, let $a_1,a_2$ be the numbers of intersection points of the curve with the
lines of the two families of the quadric. We say that it is a curve of type
$(a_1,a_2)$.  The curve $C$ has two pencils $g_{a_1}^1,g_{a_2}^1$ defined by two divisors
$\aa_1,\aa_2$. We know that
$a_1+a_2=d$ and
$g=(a_1-1)(a_2-1)-n$, where $n$ is the number of singular points with the adequate
multiplicity. 

In this way, the scroll is the image of the ruled surface
$S=\P(\Te_X(\aa_1)\oplus \Te_X(\aa_2))$ by the map defined by a linear (sub)system
$W\subset |\Te_S(1)|$. The scroll is linearly normal when $W$ is a complete linear system,
that is, when the two pencil are complete linear systems.

In the second case, let $a>1$ be the number of intersection points of the curve with a line of
the cone. If $2a<d-1$, the curve has a point of multiplicity $d-2a$ in the vertex of the
cone. Let $n$ be the number of additional singular points with the adequate multiplicity. We
know that $g=(a-1)(d-a-1)-n$. Now, the scroll has a directrix line of multiplicity $a$ and $n$ double generators.

\end{enumerate}

\subsection{Special scrolls of genus $3$ and degree $d$, $4\leq d\leq 6$.}

Let $C$ be a curve of degree $d$ and genus $3$ on the Grassmannian variety of lines in $\P^3$ which parameterizes the
generators of the scroll $R$. We will study the range $4\leq d\leq 6$:

\begin{enumerate}

\item $d=4$. Then $C$ lays on a plane $\pi$. Because $R$ is not degenerated, $\pi$ is an
$\alpha$-plane and $R$ is a cone over a linearly normal plane quartic curve of genus $3$ (the
canonical curve of genus $3$). It has speciality $4$.

\item $d=5$. Then $C$ is a plane curve. $C$ is contained in an $\alpha$-plane, so $R$ is a
cone over a nonspecial plane curve of degree $5$ and genus $3$. It has speciality $3$.

\item $d=6$. Then $C$ is on $\P^2$ or $\P^3$. $R$ is a scroll of speciality $2$.

\begin{enumerate}

\item If $C\in \P^2$ then $C$ lays on an $\alpha$-plane. $R$ is a cone over a plane curve of
degree $6$ and genus $3$. But this curve is not linearly normal so the cone is not
linearly normal.

\item If $C\in \P^3$ we have two cases:

\begin{enumerate}

\item $C$ lays on a smooth quadric. Then,

\begin{enumerate}

\item $C$ is a nonsingular curve of type $(2,4)$. The curve is hyperelliptic. It corresponds to the geometrically ruled surface $\P(\Te_X(g_2^1)\oplus \Te_X(\b))$ where $\b$ is a base point free nonspecial divisor of degree $4$.

\item $C$ is a curve of type $(3,3)$ with a singular point. The curve is not hyperelliptic. It corresponds to the geometrically ruled surface $\P(\Te_X(\k-P)\oplus \Te_X(\k-Q))$, with $P,Q\in X$, $P\neq Q$.

\end{enumerate}

\item $C$ lays on a quadric cone. Then,

\begin{enumerate}

\item $C$ meets each generator at $3$ points. It does not pass through the vertex and it has a
double point. Then $R$ is indecomposable. It has a triple line and a plane curve of degree $4$ genus $3$ meeting at a point. There is a double generator. Taking planes through this generator we obtain a $1$-dimensional family of curves of degree $4$.

\item $C$ meets each generator at $2$ points. It passes through the vertex with multiplicity
$2$ and it's smooth out of the vertex. Then $R$ is indecomposable. It has a double line and a plane curve of degree $5$ and genus $3$ meeting at a point. This point is a singular point of the plane curve with multiplicity
$3$. The directrix line is a double generator too. Taking planes through each generator we obtain
a $2$-dimensional family of curves of degree $5$. 

\end {enumerate}

\end{enumerate}

\end{enumerate}

\end{enumerate}

\subsection{Special scrolls of genus $3$ and degree $7$.}

In the previous section we have completed all the cases when $d\leq 6$. Let us study  the special
scrolls of speciality $1$ and degree $7$. They are obtained by projecting $R_{\b}$ from $3$ points
$\overline{x_1},\overline{x_2},\overline{x_3}$ laying in a plane, where $deg(\b)=6$. This is equivalent to study the elementary transform $S_{\b}'$ of $S_{\b}$ at three points $x_1,x_2,x_3$.

We know the families of unisecant curves on $S_{\b}$:

- $X_0$, with $X_0^2=-2$. 

- $X_1$, with $X_1^2=2$.

- $Y_c\equiv X_1+cf$, with $Y_c^2=2(c+1)$.

Applying the properties of the elementary transform, we can see which are the curves of minimun self-intersection of $S_{\b}'$. We call $P_i=\pi(x_i)$, with $i=1,2,3$.

\begin{enumerate}

\item If $\overline x_1,\overline x_2\in \overline{X_0}$ then the line joining both points meets $\overline{X_0}$ at four points, because $\overline{X_0}$ is a plane curve. We are not projecting exactly from three points, so we does not consider this case.

\item If $\overline x_1\in \overline{X_0}$, and $\overline{x_2},\overline{x_3}\not\in \overline{X_0}$ then the unique minimum self-intersection curve of $S_{\b}'$ is $X_0'$, with $X_0'^2=-1$. If $\overline{x_2},\overline{x_3}\in \overline{X_1}$ then $X_1'^2=1$ and $S_{\b}'$ is the decomposable ruled surface $\P(\Te_X(\k-P_1)\oplus \Te_X(\b-P_2-P_3))$. In other case, $S_{\b}'$ is an indecomposable ruled surface with invariant $e=1$.

\item If $\overline{x_1},\overline{x_2},\overline{x_3}\not\in \overline{X_0}$:

\begin{enumerate}

\item If $\overline{x_1},\overline{x_2},\overline{x_3}\in \overline{X_1}$, the the plane containing the three points meets $\overline{X_1}$ at six points. But we must project exactly from three points. We discard this case.

\item If $\overline{x_1}\not\in \overline{X_1}$ then one of the minimum self-intersection curve of $S_{\b}'$ is $X_0'$, with $X_0'^2=1$. Thus, $S_{\b}'$ is an indecomposable scroll with invariant $e=-1$. It has one or two minimum self-intersection curve depending on the position of the three points of projection. 

\end{enumerate}

\end{enumerate}

Note, that when $X$ is hyperelliptic the curve $\overline{X_0}$ is a double conic. Thus, in this case, the number of points of projection in $\overline{X_0}$ must be even.

\subsection{Classification of special scrolls of genus $3$.}

From the previous discussions  we obtain the following classification:

\begin{teo}

Let $R\subset \P^3$ be a linearly normal special scroll of genus $3$, degree $d$ and speciality $i$,
over a nonhyperelliptic curve. Then $R$ is isomorphic to one of the following models:

\footnotesize
$$
\begin{array}{|c|c|c|c|c|c|}
\hline
{d}&{i}&{\E_0}&{e}&{|H|}&{
\begin{array}{c}
{\hbox{Curves of minimun degree}}\\
{nC_k^*; \hbox{$k=$degree; $^*$=special}}\\
{\hbox{$n=$multiplicity of $C_k$ in $R$}}\\
\end{array}}\\
\hline
\hline
{4}&{4}&{\begin{array}{l}
{\Te_X\oplus\Te_X(-\k)}\\
\end{array}}&{4}&{
\begin{array}{c}
{|X_0+\k f|}\\
{\hbox{(cone)}}\\
\end{array}}&{C_4^*\subset \P^2\quad (\infty^3)}\\
\hline
\hline
{5}&{3}&{\begin{array}{l}
{\Te_X\oplus\Te_X(-\b)}\\
{deg(\b)=5}\\
{\b\hbox{ base-point-free}}\\
\end{array}}&
{5}&{
\begin{array}{c}
{|X_0+\b f|}\\
{\hbox{(cone)}}\\
\end{array}}&{C_5\subset \P^2\quad (\infty^3)}\\
\hline
\hline
{6}&{2}&{\begin{array}{l}
{\Te_X\oplus\Te_X(Q-P)}\\
{P,Q\in X}\\
{P\neq Q}\\
\end{array}}&{0}&{|X_0+(\k-P) f|}&{
\begin{array}{l}
{\phi(X_0)=3C_1^*\subset \P^1}\\
{\phi(X_1)=3C_1^*\subset \P^1}\\
\end{array}}\\
\hline
{6}&{2}&{\begin{array}{l}
{\hbox{indecomposable}}\\
{\e\sim 0}\\
\end{array}}&{0}&{|X_0+(\k-P) f|}&{
\begin{array}{l}
{\phi(X_0)=3C_1^*\subset \P^1}\\
{C_4^*\subset \P^2\quad (\infty^1)}\\
\end{array}}\\
\hline
\hline
{7}&{1}&{\begin{array}{l}
{\hbox{indecomposable}}\\
\end{array}}&{-1}&{|X_0+(\k-\e) f|}&{
\begin{array}{l}
{\phi(X_0)=C_4^*\subset \P^2}\\
{C_5\subset \P^2\quad (\infty^1)}\\
\end{array}}\\
\hline
{7}&{1}&{\begin{array}{l}
{\hbox{indecomposable}}\\
\end{array}}&{-1}&{|X_0+(\k-\e) f|}&{
\begin{array}{l}
{\phi(X_0)=C_4^*\subset \P^2}\\
{4C_1\subset \P^1}\\
\end{array}}\\
\hline
{7}&{1}&{\begin{array}{l}
{\Te_X\oplus \Te_X(\k-\b-P)}\\
{deg(\b)=4}\\
{\b\hbox{ nonspecial}}\\
{\b\hbox{ base-point-free}}\\
{P\in X}\\
\end{array}}&{1}&{|X_0+\b f|}&{
\begin{array}{l}
{\phi(X_0)=3C_1^*\subset \P^1}\\
{\phi(X_1)=4C_1\subset \P^1}\\
\end{array}}\\
\hline
{7}&{1}&{\begin{array}{l}
{\hbox{indecomposable}}\\
\end{array}}&{1}&{|X_0+(\k-P-\e) f|}&{
\begin{array}{l}
{\phi(X_0)=3C_1^*\subset \P^1}\\
{C_5\subset \P^2\quad (\infty^1)}\\
\end{array}}\\
\hline
\end{array}
$$
\normalsize
where $\phi$ denotes the map defined by the linear system $|H|$:
$$
\phi:\P(\E_0)\lrw R\subset \P^3.
$$

\end{teo}

\begin{teo}

Let $R\subset \P^3$ be a linearly normal special scroll of genus $3$, degree $d$ and speciality $i$,
over a hyperelliptic curve. Then $R$ is isomorphic to one of the following models:

\footnotesize
$$
\begin{array}{|c|c|c|c|c|c|}
\hline
{d}&{i}&{\E_0}&{e}&{|H|}&{
\begin{array}{c}
{\hbox{Curves of minimun degree}}\\
{nC_k^*; \hbox{$k=$degree; $^*$=special}}\\
{\hbox{$n=$multiplicity of $C_k$ in $R$}}\\
\end{array}}\\
\hline
\hline
{5}&{3}&{\begin{array}{l}
{\Te_X\oplus\Te_X(-\b)}\\
{deg(\b)=5}\\
{\b\hbox{ base-point-free}}\\
\end{array}}&
{5}&{
\begin{array}{c}
{|X_0+\b f|}\\
{\hbox{(cone)}}\\
\end{array}}&{C_5\subset \P^2\quad (\infty^3)}\\
\hline
\hline
{6}&{2}&{\begin{array}{l}
{\Te_X\oplus\Te_X(g_2^1-\b)}\\
{deg(\b)=4}\\
{\b\hbox{ base-point-free}}\\
{\b\hbox{ nonspecial}}\\
\end{array}}&{2}&{|X_0+\b f|}&{
\begin{array}{l}
{\phi(X_0)=2C_1^*\subset \P^1}\\
{\phi(X_1)=4C_1\subset \P^1}\\
\end{array}}\\
\hline
{6}&{2}&{\begin{array}{l}
{\hbox{indecomposable}}\\
{\e\sim 0}\\
\end{array}}&{2}&{|X_0+(g_2^1-\e) f|}&{
\begin{array}{l}
{\phi(X_0)=2C_1^*\subset \P^1}\\
{C_5\subset \P^2\quad (\infty^2)}\\
\end{array}}\\
\hline
\hline
{7}&{1}&{\begin{array}{l}
{\hbox{indecomposable}}\\
\end{array}}&{-1}&{|X_0+(\k-\e) f|}&{
\begin{array}{l}
{\phi(X_0)=2C_2^*\subset \P^2}\\
{C_5\subset \P^2\quad (\infty^1)}\\
\end{array}}\\
\hline
{7}&{1}&{\begin{array}{l}
{\hbox{indecomposable}}\\
\end{array}}&{-1}&{|X_0+(\k-\e) f|}&{
\begin{array}{l}
{\phi(X_0)=2C_2^*\subset \P^2}\\
{4C_1\subset \P^1}\\
\end{array}}\\
\hline
\end{array}
$$
\normalsize
where $\phi$ denotes the map defined by the linear system $|H|$:
$$
\phi:\P(\E_0)\lrw R\subset \P^3.
$$

\end{teo}

\end{document}